\documentclass{article}
\usepackage{mathptmx} 

\newcommand\lpp\measuredangle 

\usepackage{url} 
\usepackage{amssymb} 
\usepackage{covington} 
\newcommand\imp{\mathbin{\rightarrow}} 
\newcommand\ttt{\mathbf{t}}
\newcommand\eee{\mathbf{e}}

\newcommand\donotcirculate[1]
{} 

\renewcommand\paragraph[1]{\medskip\par\noindent\textbf{#1}\quad} 

\title{Specimens:\quad ``most of" generic NPs\\ in a contextually flexible type theory} 
\author{Christian Retor\'e (LaBRI \& INRIA, Universit\'e de Bordeaux)}
\date{Third Genius Conference, Paris, December 5-6 2011 \url{http://geniusconference.org/}} 

\begin{document} 

\maketitle 

\paragraph{Overview} This paper proposes to compute the meanings associated to sentences with generic NPs  
corresponding to 
the \emph{most of} generalized quantifier. We call these generics \emph{specimens} 
and they resemble stereotypes or prototypes in  lexical semantics.  
The \emph{meanings} are viewed as logical formulae that can be thereafter interpreted in your favorite models. 

We rather depart from the dominant Fregean single untyped  universe and go for type theory 
with hints from Hilbert $\epsilon$ calculus \cite{HilbertGrundlagen1922,espilonURL} 
and from medieval philosophy see e.g. \cite{libera1993philosophie}. 
Our type theoretic analysis bears some resemblance with on going work 
in lexical semantics. \cite{asher-webofwords,BMRjolli} 

Our model also applies to classical examples involving a class (or a generic element of this class) 
which is provided by the context. 
An outcome of this study is that, in the minimalism-contextualism debate, see e.g. \cite{SJCcerisy}, 
if one adopts a type theoretical view, terms encode the purely semantic meaning component while their typing is 
 pragmatically determined. 

\paragraph{Terms for universal and specimen generics} 
Here are two examples from the web. The first one involves a universal generic element 
and the second one  a \emph{specimen}:  

\noindent (1) The AKC notes that any dog may bite [...] \label{dog} 

\noindent (2) The Brits love France. \label{brits} 

As Hilbert  calculus shows, quantifiers, classes  and generic elements are closely related. 
He introduced $\tau x.\ A$, an element such that 
$\forall x.\ A(x)$ is equivalent to $A(\tau x.\ A(x))$: 
$A$ holds for every object if and only if it holds for this element $\tau x.\ A(x)$,
i.e. it is the universal generic associated with $A$. 
\footnote{There is the dual existential generic $\epsilon x. A$ which satisfies  $A(\epsilon x. A(x))\equiv A(\tau x. \lnot A(x))$).}  
This view is rather confidential.  
Some exceptions are the work on definite NPs with $\iota$ choice function (in particular by von Heusinger see e.g. \cite{Heusinger2007}) 
and \cite{AbrusciRetoreCLMPS}  on generalized quantifiers.  

Here we suggest to associate to any property $A$ its specimen written $\lpp x.\ A$.  
Intuitively, it enjoys all the properties that are true of \emph{most of} $A$.
Although this paper remains on the ``syntactic side of semantics", inhabited with logical formulae, 
let us give a few hints on 
how to interpret specimens and their possible reference in models. 

Properties of $\lpp x.\ A$ are the ones that are true of \emph{most of} the $A$ 
\footnote{Observe that we do not fix a 
precise ratio much larger than a half. Indeed, 
 it is a vague quantifier. However \emph{most of} as opposed to what is commonly said, is not a matter of cardinality but of measure:
for instance in maths books it is said that most of number are not primes meaning $\lim_{n\rightarrow +\infty} \pi(n)\rightarrow 0$.} 
there can be no contradiction, since when $P$ holds of the specimen, $\lnot P$ does not.
-- in Hilbertian terms it is more like a $\tau$ than like an $\epsilon$. 
For scalar values we prefer to have relations rather than functions: indeed the specimen has not, for instance,  a single height but the relation  $height(spec,x)$
is true whenever $x$ is in some interval (think of baby height charts). 

As far as proofs are concerned, we know 
some situations which enables to assert that $P$ holds of $\lpp x.\ A$:  when the universal quantifier holds,
when all the \emph{most of} properties are true of it,... But, as expected, we do not know any complete set of rules. 
We also know it can be refuted when there are only a small minority  of $A$ enjoys $P$, 
or when there is another property $Q$ disjoint from $P$ and true of $\lpp x.\ A$. 

The specimen can be foreseen in ancient and medieval logic: 
the predication on object as member of some class, 
formal generic elements with a given  ontological class, 
essential and accidental properties ... 
In particular,  predicates that apply to several 
unrelated classes were distinguished from ``homogenous" predicates that apply to (the generic element of) a class, 
-- as in Abu'l Barakat  or Dun Scott.  \cite{libera1993philosophie}

\paragraph{A flexible typed calculus with a specimen operator}  
As in Montague semantics we assume that a lexicon associates 
typed $\lambda$-terms with each word, and we start from  
a syntactic analysis (saying what applies to what). 
The logical formula depicting the meaning is obtained by applying lexical $\lambda$-terms one to another, according to the syntax. 
In addition to this montagovian term depicting argumental structure  each word is also provided 
with a finite number of $\lambda$-terms which are optionally used to convert the type when needed.
For instance the lexicon provide for some human entries a term which convert them into vehicles when needed, 
e.g  if a VP like \emph{``is parked up for the night''}  is applied to \emph{``Nic''}.

Instead of simply typed $\lambda$-calculus we use second order $\lambda$-calculus, namely Girard system F (1971), see e.g.  \cite{GLT88}. 
Base types are constant types (the usual ones of $TY_n$, $\ttt$, $\eee_i$, lots of entity types), or variable types, $\alpha$, $\beta$, ...
When $T_1$ and $T_2$ are types, so is $T_1\imp T_2$ and when $T$ 
is a type and $\alpha$ is a type variable,  $\Pi\alpha.\ T$ is a type as well -- $\alpha$ usually appears in $T$ 
but not necessarily. 

As opposed to other type theories e.g. (I)TT, the system is conceptually and formally 
extremely simple, quite powerful,... and paradox free.  

\noindent 
\begin{minipage}{\textwidth}
Term building operations include  the ones of simply typed $\lambda$-calculus: 

\noindent \fbox{$vc$}  Constants (resp. variables) of a given type $T$ are terms: $c:T$ (resp. $x:T$).

\noindent \fbox{e$\lambda$}  If  $u$ is a term of type $T_1\imp T_2$ and $v$ is a term of type $T_1$, then $u(v)$ is a term of type $T_2$. 

\noindent \fbox{i$\lambda$}  If $u$ is a term of type $T_2$ and $x$ a variable of type $T_1$, then $\lambda x.\ u$ is a term of type $T_1\imp T_2$. 

\noindent These usual operations are completed by quite similar operations handling quantification over all types: 

\noindent \fbox{$e\Lambda$} If $u$ is a term of type $\Pi\alpha.\ U$ and $T$ is a type, then $u\{T\}$ is a (specialized)  term of type $U[\alpha:=T]$ 

\noindent \fbox{$i\Lambda$}  If $u$ is a term of type $T$ and if there is no occurrence of the type variable $\alpha$ 
in the type of any free variable ($u$  works uniformly for every type $\alpha$), then $\Lambda \alpha.u$ is a term of type $\Pi\alpha.\ T$ 
(that's the universal view of $u$).  
\end{minipage}

Remember usual beta-reduction is $(\lambda x^T.\ u) t^T\rightsquigarrow u[x:=t]$ Here, beta-reduction for types and $\Lambda$ works just the same: $(\Lambda \alpha.\ u)\{T\} \rightsquigarrow u[\alpha:=T]$. 

\indent In F, instead of having a constant $\forall_\alpha$ of type $(\alpha\imp \ttt)\imp\ttt$ for every type $\alpha$ over which we would like to quantify 
we shall have one constant  $\forall$ of type $\Pi\alpha.\ (\alpha\imp \ttt)\imp\ttt$ that 
\donotcirculate{5ex}will    
 be applied to $T$ to obtain the quantifier over the type $T$:
 
$$\forall\{human\}(\lambda x^{human}. mortal^{human\imp \ttt}(x)$$

We introduce a constant  $\lpp$ of type $\Pi \alpha.\ \alpha$  mapping each property to its specimen.  
When applied to a type $T$, this constant $\lpp$ yields the element $\lpp\{T\}$  of type $T$ which is assumed to be the specimen of $T$
($\lpp\{T\}$ is the F term for $\lpp x. T$ when types and properties are identified):
 it is to be interpreted consequently when interpreting the resulting formula. 
 We could also use  the type raised version, mapping each property $A$ to the average element of type $A$
 as some did for the choice function.

\paragraph{Computing the readings: semantic terms and contextual typing} 
It is easily seen that our model 
will provide the right formula for the example (2): 

$love(\lpp\{brits\},France)$

It resembles the $\iota$ choice function, apart that it selects an element  about 
which we can assert properties but which does not exists \textsl{stricto sensu},  
as   medieval universals, Hilbert's $\tau x.\ A$, etc.  

We actually started our reflexion 
on such generics from  classical examples in the minimalism-contextualism debate. 
Such statements can be both true and false depending on the 
class in which the object is considered, which is provided by the context. 

\noindent (3)  Carlotta is tall. 

If Carlotta is a two year old girl it can be both true ("My daughter is tall and thin for a 2 year old.")
and false ("My two-year-old can't get his own cup [...] because he 
can't reach, [...]") depending on her class -- her type in our type theoretic framework. 

We noticed that the specimen notion together with the flexibility of F typing succeed to capture this phenomenon. 
Many of optional $\lambda$-terms encode the ontological relations and in the case of a two-year old girl like Carlotta, she can be viewed as  a child, 
and also as a female human being, as a human being etc. 

Here are the constants and the useful lexicon entries: 

\noindent $height:\Pi \alpha.\ (\alpha \imp \mathit{float} \imp \ttt)$ 

\noindent $<: \mathit{float} \imp \mathit{float} \imp \ttt$

\noindent \textbf{Carlotta} $Carlotta: 2yoGirl$  (constant) 

\noindent \qquad $\mathsf{h}: 2yoGirl \imp human$  (optional $\lambda$-term)

\noindent \textbf{tall} $\Lambda \alpha. \lambda x^\alpha \forall\{\mathit{float}\} \lambda h_s^{\mathit{float}} \forall\{\mathit{float}\}  \lambda h^{\mathit{float}} 
\newline\hspace*{3ex}
height\{\alpha\}(\lpp\{\alpha\},h_s)\land height\{\alpha\}(x,h)\Rightarrow h_s < h$

\noindent\hspace*{3ex}type of tall: $\Pi \alpha. \alpha\imp\ttt$ 

The constant $height$ is a relation between members of a type and numbers ($\mathit{float}$) which are compared with $<$. 
The entry for \emph{tall} applies to any type $T$ (second order is quite important here as well) and to a term $u$ of type $T$.
It says that the object $u$ is taller than any possible height of the specimen of this class $T$. 

If we do not use any optional $\lambda$-term, we apply \emph{tall} to the type $2yoGirl$, and to the constant $Carlotta^{2yoGirl}$
 we get the reading where Carlotta is taller than the maximal height of the $2yoGirl$
specimen  (think again of baby height charts). This is likely to be interpreted as true. 
 
But if we apply \emph{tall} to the \emph{human} type, we cannot apply the result to the constant $Carlotta^{2yoGirl}$.
But we can apply the $\mathsf{h}: 2yoGirl \imp human$  (optional $\lambda$-term) to the constant $Carlotta^{2yoGirl}$ and proceed: 
using the type $human$ since $\mathsf{h}(Carlotta)$ is of type $human$. 
We thus obtain the formula meaning that Carlotta is tall as a human being, which is unlikely to be interpreted as true. 

The semantic machinery produces every possible reading 
and the context intervenes as a preference for some optional transformation(s). 
It should be discussed whether there is one or several natural types for an object. 
Our model can handle any solution: a single natural type,  several privileged types,...
--- quite often, such ontological or metaphysical questions 
spontaneously pop up when dealing with the organization of the concepts in the lexicon. 
 
\paragraph{Conclusion} We applied the F typed $\lambda$-calculus to derive 
semantic readings in the presence of \emph{``most of"} generics, that we call \emph{specimens}. 
Our treatment also helps to determine the border between 
semantics and pragmatics: 
the term calculus models the semantics, while the typing flexibility of F represents the possible context adaptation.  

\paragraph{\itshape Thanks} \emph{This work owes a lot to Sarah-Jane Conrad (\url{meaning.ch}, Sprachphilosophie, Universit\"at Bern). 
Indeed, her talk and our discussions initiated in Cerisy on the debate between contextualism 
and semantic minimalism, lead us to a new connection between logical semantics and type theory, here applied to generic elements.}  
\bibliography{bigbiblio} 
\end{document}